\documentclass[a4paper, 12pt]{article}%
\usepackage{amsmath}
\usepackage{graphicx}%
\usepackage{amsfonts}%
\usepackage{amssymb}

\begin{document}

\begin{center}
{\Huge Joinings of W*-dynamical systems}

\bigskip

Rocco Duvenhage

\bigskip

\textit{Department of Mathematics and Applied Mathematics}

\textit{University of Pretoria, 0002 Pretoria, South Africa}

\bigskip

rocco.duvenhage@up.ac.za

\bigskip

2007-11-15
\end{center}

\bigskip

\noindent\textbf{Abstract}

We study the notion of joinings of W*-dynamical systems, building on ideas
from measure theoretic ergodic theory. In particular we prove sufficient and
necessary conditions for ergodicity in terms of joinings, and also briefly
look at conditional expectation operators associated with joinings.

\bigskip


\section{Introduction}

The study of joinings (and disjointness) of measure theoretic dynamical
systems was initiated by Furstenberg \cite{F67} in 1967, and Rudolph
\cite{Rud} in 1979. Joinings have since become a useful tool in ergodic
theory. More recent treatments of joinings including further developments and
some applications can be found in Glasner's book \cite{G}, Rudolph's book
\cite{Rud90}, the review \cite{dlR}, and the paper \cite{LPT}.

In this paper we study joinings of W*-dynamical systems. We will refer to
W*-dynamical systems simply as ``dynamical systems''; see Section 2 for the
precise definition that we'll use. In these dynamical systems one works on a
von Neumann algebra rather than a measurable space, and with a state instead
of a measure. The von Neumann algebra is a noncommutative generalization of
the abelian algebra $L^{\infty}$ in the measure theoretic case. Such
noncommutative dynamical systems have of course been studied extensively, and
is for example a suitable framework for the mathematical study of quantum
physics. Roughly a joining of two dynamical systems is a generalization of the
usual product of the systems, but where the product state is replaced with a
state which can in principle take into account an ``overlap'' or ``common
part'' of the two systems. It is such a state that we will refer to as a
``joining'' of the two dynamical systems. We focus mainly on sufficient and
necessary conditions for ergodicity in terms of joinings (see Section 3), but
also consider some basic aspects of conditional expectation operators
associated with joinings.

The same idea has been used by Sauvageot and Thouvenot \cite{ST} to study
entropy in noncommutative dynamical systems. However they consider joinings
where one of the two systems is classical (i.e. its algebra is abelian). A
useful reference regarding this approach to entropy is \cite[Chapter 5]{St}.
Although we will not apply our results to entropy in this paper, this joining
approach to entropy (together with the work on joinings in classical ergodic
theory) does suggest that a general study of joinings of noncommutative
dynamical systems will have uses beyond just ergodicity. We also note that in
the literature on entropy the term ``stationary coupling'' is used, rather
than ``joining'', however the latter is standard in measure theoretic ergodic
theory, more succinct, and appears to be older, so we will continue to use it.

For the most part we consider general group actions, but a necessary condition
for ergodicity is only proved in the case of amenable countable discrete
groups. For the proof of the sufficient condition it is useful to work in
terms of a ``factor'' (essentially a subsystem) of the dynamical system, and
the definition of a factor is given in Section 3. To avoid confusion, note
that in this context the term does not refer to a von Neumann algebra $A$
which is a factor (i.e. $A\cap A^{\prime}=\mathbb{C}1$).

Our von Neumann algebras always contain the identity operator on the
underlying Hilbert space, and we will denote it by $1_{A}$ for a von Neumann
algebra $A$, or sometimes just $1$. The identity map $A\rightarrow A$ will be
denoted by id$_{A}$ or simply id, while the group of all $\ast$-automorphisms
of $A$ will be denoted by Aut$(A)$. We only consider dynamical systems on
$\sigma$-finite von Neumann algebras, since we will be using Tomita-Takesaki
theory in Section 3. Remember that a von Neumann algebra is $\sigma$-finite if
and only if it has a faithful normal state. We will denote the algebraic
tensor product of two von Neumann algebras $A$ and $B$ by $A\odot B$, which is
a unital $\ast$-algebra. For simplicity we will only consider algebraic tensor
products in this paper. The von Neumann algebra of bounded linear operators
$H\rightarrow H$ on a Hilbert space $H$ will be denoted by $B(H)$, and the
commutant of any $S\subset B(H)$ is denoted by $S^{\prime}$. Our main
reference for von Neumann algebras and Tomita-Takesaki theory is \cite{BR}.

Throughout this paper $G$ is\ an arbitrary but fixed group, except in Theorem
3.7 where we specialize.

\section{Joinings}

This section is devoted mainly to the basic definitions, and includes a
characterization of joinings in terms of conditional expectation operators.

\bigskip

\textbf{Definition 2.1.} A \textit{dynamical system} $\mathbf{A}=\left(
A,\mu,\alpha\right)  $ consists of a faithful normal state $\mu$ on a $\sigma
$-finite von Neumann algebra $A$, and a representation $\alpha:G\rightarrow$
Aut$(A):g\mapsto\alpha_{g}$ of $G$ as $\ast$-automorphisms of $A$, such that
$\mu\circ\alpha_{g}=\mu$ for all $g$. We will call $\mathbf{A}$
\textit{trivial} if $A=\mathbb{C}1_{A}$. We will call $\mathbf{A}$ an
\textit{identity system} if $\alpha_{g}=$ id$_{A}$ for all $g$.

\bigskip

In the remainder of this section and the next, the symbols $\mathbf{A}$,
$\mathbf{B}$ and $\mathbf{F}$ will denote dynamical systems $\left(
A,\mu,\alpha\right)  $, $\left(  B,\nu,\beta\right)  $ and $\left(
F,\lambda,\varphi\right)  $ respectively, and keep in mind that they\ all make
use of actions of the same group $G$.

\bigskip

\textbf{Definition 2.2.} A \textit{joining} of $\mathbf{A}$ and $\mathbf{B}$
is a state $\omega$ on $A\odot B$ such that%
\begin{align*}
\omega\left(  a\otimes1_{B}\right)   &  =\mu(a)\\
\omega\left(  1_{A}\otimes b\right)   &  =\nu(b)
\end{align*}
and
\[
\omega\circ\left(  \alpha_{g}\otimes\beta_{g}\right)  =\omega
\]
for all $a\in A$, $b\in B$ and $g\in G$. The set of all joinings of
$\mathbf{A}$ and $\mathbf{B}$ is denoted by $J\left(  \mathbf{A}%
,\mathbf{B}\right)  $. Note that $\mu\otimes\nu\in J\left(  \mathbf{A}%
,\mathbf{B}\right)  $. We call $\mathbf{A}$ \textit{disjoint} from
$\mathbf{B}$ when $J\left(  \mathbf{A},\mathbf{B}\right)  =\left\{  \mu
\otimes\nu\right\}  $.

\bigskip

As part of the proof of Theorem 3.3 in Section 3, we will construct a joining
other than $\mu\otimes\nu$ in the special case where $\mathbf{B}$ is obtained
in a certain way from a ``factor'' of $\mathbf{A}$.

We are now going to study the conditional expectation operator associated with
certain states on $A\odot B$.

\bigskip

\textbf{Construction 2.3.} Let $\omega$ be any state on $A\odot B$ such that
$\omega\left(  a\otimes1_{B}\right)  =\mu(a)$ and $\omega\left(  1_{A}\otimes
b\right)  =\nu(b)$.

Consider the GNS construction $\left(  H_{\omega},\gamma_{\omega}\right)  $
for $\left(  A\odot B,\omega\right)  $, by which we mean $H_{\omega}$ is a
Hilbert space and
\[
\gamma_{\omega}:A\odot B\rightarrow H_{\omega}%
\]
a linear operator such $\gamma_{\omega}\left(  A\odot B\right)  $ is dense in
$H_{\omega}$ and $\left\langle \gamma_{\omega}(s),\gamma_{\omega
}(t)\right\rangle =\omega\left(  s^{\ast}t\right)  $ for all $s,t\in A\odot
B$. Then $\Omega_{\omega}:=\gamma_{\omega}\left(  1_{A}\otimes1_{B}\right)  $
is the corresponding cyclic vector.

Define $\iota_{A}:A\rightarrow A\odot B:a\mapsto a\otimes1_{B}$ and $\iota
_{B}:B\rightarrow A\odot B:b\mapsto1_{A}\otimes b$ and let $H_{\mu}$ and
$H_{\nu}$ be the closures in $H_{\omega}$ of $\gamma_{\omega}\circ\iota
_{A}(A)$ and $\gamma_{\omega}\circ\iota_{B}(B)$ respectively. Setting
\[
\gamma_{\mu}:=\gamma_{\omega}\circ\iota_{A}:A\rightarrow H_{\mu}%
\]
and
\[
\gamma_{\nu}:=\gamma_{\omega}\circ\iota_{B}:B\rightarrow H_{\nu}%
\]
we have $\gamma_{\mu}(A)$ and $\gamma_{\nu}(B)$ dense in $H_{\mu}$ and
$H_{\nu}$ respectively, and $\left\langle \gamma_{\mu}(a),\gamma_{\mu
}(a^{\prime})\right\rangle =\omega\left(  \left(  a\otimes1_{B}\right)
^{\ast}\left(  a^{\prime}\otimes1_{B}\right)  \right)  =\mu\left(  a^{\ast
}a^{\prime}\right)  $ for all $a,a^{\prime}\in A$, and similarly $\left\langle
\gamma_{\nu}(b),\gamma_{\nu}(b^{\prime})\right\rangle =\nu\left(  b^{\ast
}b^{\prime}\right)  $. Hence $\left(  H_{\mu},\gamma_{\mu}\right)  $ and
$\left(  H_{\nu},\gamma_{\nu}\right)  $ are the GNS constructions for $\left(
A,\mu\right)  $ and $\left(  B,\nu\right)  $ respectively, and they both have
the cyclic vector $\Omega_{\mu}:=\gamma_{\mu}\left(  1_{A}\right)
=\Omega_{\omega}=\gamma_{\nu}\left(  1_{B}\right)  =:\Omega_{\nu}$.

Let $P$ be the projection of $H_{\omega}$ onto the subspace $H_{\nu}$ and then
set
\[
P_{\omega}:=P|_{H_{\mu}}:H_{\mu}\rightarrow H_{\nu}%
\]
which is called the \textit{conditional expectation operator} associated with
$\omega$. It is the unique mapping $H_{\mu}\rightarrow H_{\nu}$ satisfying
\[
\left\langle P_{\omega}x,y\right\rangle =\left\langle x,y\right\rangle
\]
for all $x\in H_{\mu}$ and $y\in H_{\nu}$. The space of fixed points of
$P_{\omega}$ is clearly $H_{\mu}\cap H_{\nu}$. In particular $P_{\omega}%
\Omega_{\omega}=\Omega_{\omega}$.

Since $\mu\circ\alpha_{g}=\mu$ and $\nu\circ\beta_{g}=\nu$, we obtain well
defined and unique linear operators $U_{g}:H_{\mu}\rightarrow H_{\mu}$ and
$V_{g}:H_{\nu}\rightarrow H_{\nu}$ from $U_{g}\gamma_{\mu}(a):=\gamma_{\mu
}\left(  \alpha_{g}(a)\right)  $ and $V_{g}\gamma_{\nu}(b):=\gamma_{\nu
}\left(  \beta_{g}(b)\right)  $. For the same reason $U_{g}$ and $V_{g}$ are
isometries. Furthermore by uniqueness, $g\mapsto U_{g}$ and $g\mapsto V_{g}$
are representations of $G$, since $\alpha$ and $\beta$ are. In particular
$U_{g}$ and $V_{g}$ are invertible, and hence unitary.

Note that this whole construction goes through even if we only assume that $A
$ and $B$ are unital $\ast$-algebras rather than von Neumann algebras.

If we furthermore assume $\omega\in J\left(  \mathbf{A},\mathbf{B}\right)  $,
which means we additionally have $\omega\circ\left(  \alpha_{g}\otimes
\beta_{g}\right)  =\omega$, then in the same way we obtain a unitary
representation $g\mapsto W_{g}$ of $G$ on $H_{\omega}$ such that $W_{g}%
\gamma_{\omega}(t)=\gamma_{\omega}\left(  \alpha_{g}\otimes\beta
_{g}(t)\right)  $ for all $t\in A\odot B$. Note that $W_{g}|_{H_{\mu}}=U_{g}$
and $W_{g}|_{H_{\nu}}=V_{g}$. $\square$

\bigskip

\textbf{Proposition 2.4.}\textit{\ Let }$\omega$\textit{\ be a state on
}$A\odot B$\textit{\ such that }$\omega\left(  a\otimes1_{B}\right)  =\mu
(a)$\textit{\ and }$\omega\left(  1_{A}\otimes b\right)  =\nu(b)$\textit{\ for
all }$a\in A$\textit{\ and }$b\in B$\textit{. Then }$\omega\in J\left(
\mathbf{A},\mathbf{B}\right)  $\textit{\ if and only if }%
\[
P_{\omega}U_{g}=V_{g}P_{\omega}%
\]
\textit{for all }$g$\textit{, in terms of Construction 2.3.}

\bigskip

\textbf{Proof.} Assuming $\omega\in J\left(  \mathbf{A},\mathbf{B}\right)  $,
then by Construction 2.3
\[
\left\langle V_{g}^{\ast}P_{\omega}U_{g}x,y\right\rangle =\left\langle
P_{\omega}W_{g}x,W_{g}y\right\rangle =\left\langle W_{g}x,W_{g}y\right\rangle
=\left\langle x,y\right\rangle =\left\langle P_{\omega}x,y\right\rangle
\]
for all $x\in H_{\mu}$ and $y\in H_{\nu}$, hence $V_{g}^{\ast}P_{\omega}%
U_{g}=P_{\omega}$. Conversely, if $P_{\omega}U_{g}=V_{g}P_{\omega}$, then for
all $a\in A$ and $b\in B$
\begin{align*}
\omega\left(  \alpha_{g}\otimes\beta_{g}(a\otimes b)\right)   &
=\omega\left(  \alpha_{g}\otimes\beta_{g}(a^{\ast}\otimes1_{B})^{\ast}%
\alpha_{g}\otimes\beta_{g}(1_{A}\otimes b)\right) \\
&  =\left\langle U_{g}\gamma_{\mu}\left(  a^{\ast}\right)  ,V_{g}\gamma_{\nu
}\left(  b\right)  \right\rangle \\
&  =\left\langle P_{\omega}U_{g}\gamma_{\mu}\left(  a^{\ast}\right)
,V_{g}\gamma_{\nu}\left(  b\right)  \right\rangle \\
&  =\left\langle P_{\omega}\gamma_{\mu}\left(  a^{\ast}\right)  ,\gamma_{\nu
}\left(  b\right)  \right\rangle \\
&  =\left\langle \gamma_{\mu}\left(  a^{\ast}\right)  ,\gamma_{\nu}\left(
b\right)  \right\rangle \\
&  =\omega\left(  a\otimes b\right)
\end{align*}
so $\omega\circ\left(  \alpha_{g}\otimes\beta_{g}\right)  =\omega$ by
linearity. $\square$

\section{Ergodicity}

We now turn to ergodicity, in particular proving sufficient and necessary
conditions for ergodicity in terms of joinings. As part of the proof of
sufficiency (Theorem 3.3) we construct a special joining in terms of a factor
of a dynamical system. The commutant of the algebra, and the modular
conjugation operator from Tomita-Takesaki theory play a central role in this construction.

\bigskip

\textbf{Definition 3.1.} A dynamical system $\mathbf{A}$ is called
\textit{ergodic} if its \textit{fixed point algebra}%
\[
A_{\alpha}:=\left\{  a\in A:\alpha_{g}(a)=a\text{ for all }g\in G\right\}
\]
is trivial, i.e. $A_{\alpha}=\mathbb{C}1_{A}$.

\bigskip

\textbf{Definition 3.2.} We call $\mathbf{F}$ a \textit{factor} of
$\mathbf{A}$ if there exists an injective unital $\ast$-homomorphism $h$ of
$F$ onto a von Neumann subalgebra of $A$ such that $\mu\circ h=\lambda$ and
$\alpha_{g}\circ h=h\circ\varphi_{g}$ for all $g\in G$. If this factor is an
identity system, then we will call it an \textit{identity factor}.

\bigskip

It is easily seen that $A_{\alpha}$ is itself a $\sigma$-finite von Neumann
algebra with $\mu|_{A_{\alpha}}$ a faithful normal state, and that
$\mathbf{A}_{\alpha}:=\left(  A_{\alpha},\mu|_{A_{\alpha}},\alpha|_{A_{\alpha
}}\right)  $ is an identity factor of $\mathbf{A}$.

\bigskip

\textbf{Theorem 3.3.}\textit{\ If }$\mathbf{A}$\textit{\ is disjoint from all
identity systems, then it is ergodic.}

\bigskip

In order to prove this theorem, we will use a special case of the following construction:

\bigskip

\textbf{Construction 3.4.} Let $\mathbf{F}$ be any factor of $\mathbf{A}$
given by the $\ast$-homomorphism $h:F\rightarrow A$ as in Definition 3.2.

Denote the cyclic representation of $\left(  A,\mu\right)  $, obtained using
the GNS construction, by $\left(  H,\pi,\Omega\right)  $. For every $g\in G$
there is a unique unitary operator $U_{g}:H\rightarrow H$ such that
$U_{g}\Omega=\Omega$ and
\[
U_{g}\pi(a)U_{g}^{\ast}=\pi\left(  \alpha_{g}(a)\right)
\]
for all $a\in A$. The uniqueness ensures that $g\mapsto U_{g}$ is a
representation of $G$.

Since $\mu$ is faithful and normal, $\Omega$ is a cyclic and separating vector
for the von Neumann algebra $M:=\pi(A)$ and $\pi:A\rightarrow M$ is a $\ast
$-isomorphism. It also follows that $\pi$ and its inverse are $\sigma$-weakly
continuous, hence $\pi\left(  h(F)\right)  $ is a von Neumann subalgebra of
$M$.

Let $J$ be the modular conjugation associated with $\left(  M,\Omega\right)  $
as obtained in Tomita-Takesaki theory. Remember that $J$ is anti-unitary,
$J^{2}=1$ (i.e. $J^{\ast}=J$) and $J\Omega=\Omega$. Define%
\[
j:B(H)\rightarrow B(H):a\mapsto Ja^{\ast}J
\]
then by Tomita-Takesaki theory
\[
j(M)=M^{\prime}%
\]
and furthermore $j$ is an anti-$\ast$-isomorphism, i.e. it is a linear
bijection such that $j(a^{\ast})=j(a)^{\ast}$ and $j(ab)=j(b)j(a)$ for all
$a,b\in B(H)$. Also, $j^{2}=$ id. From these facts it is easily seen that
\[
j(S)^{\prime}=j\left(  S^{\prime}\right)
\]
for all $S\subset B(H)$.

Set
\[
\sigma:=j\circ\pi\circ h
\]
then $\sigma(F)^{\prime\prime}=j\left(  \pi\left(  h(F)\right)  ^{\prime
\prime}\right)  =\sigma(F)$, since $\pi\left(  h(F)\right)  $ is a von Neumann
algebra, hence
\[
B:=\sigma(F)\subset M^{\prime}%
\]
is a von Neumann algebra. We can define a state $\nu$ on $B$ by%
\[
\nu(b):=\left\langle \Omega,b\Omega\right\rangle
\]
then clearly $\nu$ is $\sigma$-weakly continuous, i.e. normal. Furthermore,
$\nu$ is faithful, since $0=\nu\left(  b^{\ast}b\right)  =\left\|
b\Omega\right\|  ^{2}$ implies that $b=0$ because $\Omega$ is separating for
$M^{\prime}$. Now set
\begin{align*}
\beta_{g}(b)  &  :=j\circ\pi\circ\alpha_{g}\circ\pi^{-1}\circ j(b)\\
&  =JU_{g}JbJU_{g}^{\ast}J
\end{align*}
for all $b\in B$. Then it is clear that $\beta$ is a representation of $G$ as
$\ast$-automorphisms of $B$, and since $U_{g}^{\ast}\Omega=\Omega$, we have
\[
\nu\circ\beta_{g}(b)=\left\langle \Omega,JU_{g}Jb\Omega\right\rangle
=\left\langle U_{g}Jb\Omega,\Omega\right\rangle =\nu(b)
\]
for all $b\in B$. Therefore
\[
\mathbf{B}:=\left(  B,\nu,\beta\right)
\]
is a dynamical system.

Note that $\mathbf{B}$ is the ``mirror image'' of $\mathbf{F}$ in $M^{\prime}
$ in the sense that they can be said to be anti-isomorphic: $\sigma
:F\rightarrow B$ is an anti-$\ast$-isomorphism, since $j$ is. Furthermore
\[
\nu\circ\sigma=\lambda
\]
and
\[
\beta_{g}\circ\sigma=\sigma\circ\varphi_{g}%
\]
for all $g$.

We now construct a joining of $\mathbf{A}$ and $\mathbf{B}$. Consider the
bilinear mapping
\[
A\times B\rightarrow B(H):(a,b)\mapsto\pi(a)b
\]
and extend it to the linear mapping $\delta:A\odot B\rightarrow B(H)$, which
is a unital $\ast$-homomorphism, since $\pi(A)=M$ while $B\subset M^{\prime}$.
Thus we can define a state $\omega$ on $A\odot B$ by
\[
\omega(t):=\left\langle \Omega,\delta(t)\Omega\right\rangle
\]
for all $t\in A\odot B$. Then
\[
\omega\left(  a\otimes1_{B}\right)  =\left\langle \Omega,\pi(a)\Omega
\right\rangle =\mu(a)
\]
and
\[
\omega\left(  1_{A}\otimes b\right)  =\left\langle \Omega,\pi\left(
1_{A}\right)  b\Omega\right\rangle =\nu(b)
\]
for all $a\in A$ and $b\in B$. The theory of self-dual cones and standard
forms in Tomita-Takesaki theory provides (see \cite[Corollary 2.5.32]{BR}) a
unitary representation Aut$(M)\ni\theta\mapsto u(\theta)$ of the group
Aut$(M)$ on the Hilbert space $H$ such that (among other properties)
$u(\theta)au(\theta)^{\ast}=\theta(a)$ and $u(\theta)J=Ju(\theta)$ for all
$a\in M$ and $\theta\in$ Aut$(M)$, while $u(\theta)\Omega=\Omega$ for all
$\theta\in$ Aut$(M)$ for which $\left\langle \Omega,\theta(a)\Omega
\right\rangle =\left\langle \Omega,a\Omega\right\rangle $ for all $a\in M$.
Since $U_{g}$ is the unique unitary operator on $H$ satisfying $U_{g}%
\pi(a)U_{g}^{\ast}=\pi\left(  \alpha_{g}(a)\right)  $ and $U_{g}\Omega=\Omega
$, we must have $u\left(  \pi\circ\alpha_{g}\circ\pi^{-1}\right)  =U_{g}$ and
therefore
\[
U_{g}J=JU_{g}%
\]
for all $g$. Hence
\begin{align*}
\omega\circ\left(  \alpha_{g}\otimes\beta_{g}\right)  (a\otimes b)  &
=\left\langle \Omega,\pi\left(  \alpha_{g}(a)\right)  \beta_{g}(b)\Omega
\right\rangle \\
&  =\left\langle \Omega,U_{g}\pi(a)U_{g}^{\ast}JU_{g}Jb\Omega\right\rangle \\
&  =\left\langle U_{g}^{\ast}\Omega,\pi(a)U_{g}^{\ast}U_{g}JJb\Omega
\right\rangle \\
&  =\left\langle \Omega,\pi(a)b\Omega\right\rangle \\
&  =\omega(a\otimes b)
\end{align*}
and therefore by linearity $\omega\circ\left(  \alpha_{g}\otimes\beta
_{g}\right)  =\omega$. So $\omega$ is indeed a joining of $\mathbf{A}$ and
$\mathbf{B}$. $\square$

\bigskip

\textbf{Lemma 3.5.}\textit{\ In Construction 3.4 we have }$\omega=\mu
\otimes\nu$\textit{\ if and only if }$\mathbf{F}$\textit{\ is trivial.}

\bigskip

\textbf{Proof.} First note that $\mathbf{F}$ is trivial if and only if
$\mathbf{B}$ is. Now, if $\mathbf{B}$ is trivial, i.e. $B=\mathbb{C}1$, then
$\omega(a\otimes b)=\left\langle \Omega,\pi(a)b\Omega\right\rangle
=\left\langle \Omega,\pi(a)\Omega\right\rangle b=\mu(a)\nu(b)=\mu\otimes
\nu(a\otimes b)$, since we can view $b\in B$ as an element of $\mathbb{C}$. By
linearity it follows that $\omega=\mu\otimes\nu$.

Conversely, suppose $\omega=\mu\otimes\nu$. Then for any $a\in A$ and $b\in B
$
\begin{align*}
\left\langle \pi(a)\Omega,b\Omega\right\rangle  &  =\left\langle \Omega
,\pi(a^{\ast})b\Omega\right\rangle \\
&  =\omega\left(  a^{\ast}\otimes b\right) \\
&  =\mu(a^{\ast})\nu(b)\\
&  =\left\langle \pi(a)\Omega,\left\langle \Omega,b\Omega\right\rangle
\Omega\right\rangle
\end{align*}
but $\pi(A)\Omega$ is dense in $H$, and $\Omega$ is separating for $B$, hence
$b=\left\langle \Omega,b\Omega\right\rangle 1\in\mathbb{C}1$. $\square$

\bigskip

\textbf{Proof of Theorem 3.3.} Let $\mathbf{F=A}_{\alpha}$ in Construction
3.4, then $\mathbf{F}$ is an identity factor of $\mathbf{A}$ as mentioned
previously, and so $\mathbf{B}$ is an identity system. If $\mathbf{A}$ is not
ergodic, then by definition $\mathbf{F}$ is not trivial, hence $J\left(
\mathbf{A},\mathbf{B}\right)  \neq\left\{  \mu\otimes\nu\right\}  $ by Lemma
3.5 and Construction 3.4. This means that $\mathbf{A}$ is not disjoint from
$\mathbf{B}$. $\square$

\bigskip

Before we proceed to necessary conditions for ergodicity, which require
additional assumptions on the group and the allowed joinings, we briefly
return to the conditional expectation operator of Construction 2.3 for a
related but independent result:

\bigskip

\textbf{Propostion 3.6.}\textit{\ Let }$P_{\omega}$\textit{\ be as in
Construction 2.3, with }$\omega\in J\left(  \mathbf{A},\mathbf{B}\right)
$\textit{, and assume that }$\mathbf{A}$\textit{\ is ergodic and }$\mathbf{B}%
$\textit{\ an identity system. Then the fixed point space of }$P_{\omega}%
$\textit{\ is }$\mathbb{C}\Omega_{\omega}$\textit{.}

\bigskip

\textbf{Proof.} Since $\mathbf{A}$ is ergodic, the fixed point space of
$U_{G}$ is $\mathbb{C}\Omega_{\omega}$; see for example \cite[Theorem
4.3.20]{BR}. But $V_{g}=$ id, since $\mathbf{B}$ is an identity system, so for
any $x\in H_{\mu}\cap H_{\nu}$ one has $U_{g}x=W_{g}x=V_{g}x=x$, since
$\omega$ is a joining. Therefore $H_{\mu}\cap H_{\nu}=\mathbb{C}\Omega
_{\omega}$. $\square$

\bigskip

Thus far we haven't required joinings to be $\sigma$-weakly continuous, but
$\sigma$-weak continuity is of course a natural assumption in the von Neumann
algebra context, and in the next result we indeed need it.

\bigskip

\textbf{Theorem 3.7.}\textit{\ Let }$G$\textit{\ be amenable, countable and
discrete. Assume }$\omega\in J\left(  \mathbf{A},\mathbf{B}\right)
$\textit{\ is }$\sigma$\textit{-weakly continuous. If }$\mathbf{A}%
$\textit{\ is ergodic and }$\mathbf{B}$\textit{\ an identity system, then
}$\omega=\mu\otimes\nu$\textit{.}

\bigskip

\textbf{Proof.} We follow a standard plan from measure theoretic ergodic
theory as can be found in \cite[Proposition 2.2]{dlR}. Let $\left(
\Lambda_{n}\right)  $ be a (right) F\o lner sequence in $G$, i.e. every
$\Lambda_{n}$ is a compact (in other words, finite) subset of $G$ with
$\left|  \Lambda_{n}\right|  >0$ such that%
\[
\lim_{n\rightarrow\infty}\frac{\left|  \Lambda_{n}\bigtriangleup\left(
\Lambda_{n}g\right)  \right|  }{\left|  \Lambda_{n}\right|  }=0
\]
for all $g\in G$ (see for example \cite[Theorems 1 and 2]{E74} for the general
theory). For any $a\in A$ and $b\in B$ we then have
\begin{align*}
\omega\left(  a\otimes b\right)   &  =\frac{1}{\left|  \Lambda_{n}\right|
}\sum_{g\in\Lambda_{n}}\omega\left(  a\otimes b\right) \\
&  =\frac{1}{\left|  \Lambda_{n}\right|  }\sum_{g\in\Lambda_{n}}\omega
\circ\left(  \alpha_{g}\otimes\beta_{g}\right)  \left(  a\otimes b\right) \\
&  =\omega\left(  \left(  \frac{1}{\left|  \Lambda_{n}\right|  }\sum
_{g\in\Lambda_{n}}\alpha_{g}(a)\right)  \otimes b\right)  \text{.}%
\end{align*}
Let $\left(  H,\pi,\Omega\right)  $ be the cyclic representation of $\left(
A,\mu\right)  $ obtained from the GNS construction, and $g\mapsto U_{g}$ the
corresponding unitary representation of $G$ on $H$ obtained from $\alpha$, and
set $\gamma=\pi(\cdot)\Omega$. Consider any $c\in\pi(A)^{\prime}$, then%
\begin{align*}
\pi\left(  \frac{1}{\left|  \Lambda_{n}\right|  }\sum_{g\in\Lambda_{n}}%
\alpha_{g}(a)-\mu(a)1_{A}\right)  c\Omega &  =c\left(  \frac{1}{\left|
\Lambda_{n}\right|  }\sum_{g\in\Lambda_{n}}U_{g}\gamma(a)-\left(
\Omega\otimes\Omega\right)  \gamma(a)\right) \\
&  \rightarrow0
\end{align*}
by the mean ergodic theorem, since $\mathbf{A}$ is ergodic and hence the fixed
point space of $U_{G}$ is $\mathbb{C}\Omega$, which corresponds to the
projection $\Omega\otimes\Omega$. Since $\Omega$ is cyclic for $\pi
(A)^{\prime}$, i.e. $\pi(A)^{\prime}\Omega$ is dense in $H$, while $\pi\left(
\frac{1}{\left|  \Lambda_{n}\right|  }\sum_{g\in\Lambda_{n}}\alpha_{g}%
(a)-\mu(a)1_{A}\right)  $ is a bounded sequence, it follows that this sequence
converges strongly and hence weakly to $0$. However the weak and $\sigma$-weak
topologies are the same on bounded norm closed balls, hence the sequence
converges $\sigma$-weakly to $0$. But $\pi^{-1}$ is a $\ast$-isomorphism
between von Neumann algebras, and hence $\sigma$-weakly continuous, therefore
\[
e_{n}:=\frac{1}{\left|  \Lambda_{n}\right|  }\sum_{g\in\Lambda_{n}}\alpha
_{g}(a)-\mu(a)1_{A}%
\]
converges $\sigma$-weakly (and hence weakly) to $0$. If the Hilbert spaces on
which $A$ and $B$ are defined are denoted $H_{A}$ and $H_{B}$ respectively,
then we therefore have $\left\langle x,e_{n}y\right\rangle \rightarrow0$ for
all $x,y\in H_{A}$ hence $\left\langle x_{1}\otimes x_{2},e_{n}\otimes
b\left(  y_{1}\otimes y_{2}\right)  \right\rangle =\left\langle x_{1}%
,e_{n}y_{1}\right\rangle \left\langle x_{2},by_{2}\right\rangle \rightarrow0$
for all $x_{1},y_{1}\in H_{A}$ and $x_{2},y_{2}\in H_{B}$. Since $\left(
e_{n}\right)  $ is bounded, and the finite linear combinations of elementary
tensors are dense in $H_{A}\otimes H_{B}$, it follows that $e_{n}\otimes b$
converges weakly, and hence $\sigma$-weakly because of boundedness, to $0$.
This means $\omega\left(  e_{n}\otimes b\right)  \rightarrow0$, from which we
conclude that $\omega\left(  a\otimes b\right)  =\omega\left(  \mu
(a)1_{A}\otimes b\right)  =\mu(a)\omega\left(  1_{A}\otimes b\right)
=\mu(a)\nu(b)$, so $\omega=\mu\otimes\nu$. $\square$

\bigskip

\noindent\textbf{Acknowledgments. }I thank Conrad Beyers and Richard de Beer
for very useful conversations.


\begin{thebibliography}{9}                                                                                                %

\bibitem {BR}O. Bratteli, D. W. Robinson, Operator algebras and quantum
statistical mechanics 1, second edition, Springer-Verlag, New York, 1987.

\bibitem {dlR}T. de la Rue, An introduction to joinings in ergodic theory,
Discrete Contin. Dyn. Syst. 15 (2006), 121--142.

\bibitem {E74}W. R. Emerson, Large symmetric sets in amenable groups and the
individual ergodic theorem, Amer. J. Math. 96 (1974), 242--247.

\bibitem {F67}H. Furstenberg, Disjointness in ergodic theory, minimal sets,
and a problem in Diophantine approximation, Math. Systems Theory 1 (1967), 1--49.

\bibitem {G}E. Glasner, Ergodic theory via joinings, Mathematical Surveys and
Monographs 101, American Mathematical Society, Providence, RI, 2003.

\bibitem {LPT}M. Lema\'{n}czyk, F. Parreau, J.-P. Thouvenot, Gaussian
automorphisms whose ergodic self-joinings are Gaussian, Fund. Math. 164
(2000), 253--293.

\bibitem {Rud}D. J. Rudolph, An example of a measure preserving map with
minimal self-joinings, and applications, J. Analyse Math. 35 (1979), 97--122.

\bibitem {Rud90}D. J. Rudolph, Fundamentals of measurable dynamics. Ergodic
theory on Lebesgue spaces, Oxford Science Publications, The Clarendon Press,
Oxford University Press, New York, 1990.

\bibitem {ST}J.-L. Sauvageot, J.-P. Thouvenot, Une nouvelle d\'{e}finition de
l'entropie dynamique des syst\`{e}mes non commutatifs, Comm. Math. Phys. 145
(1992), 411--423.

\bibitem {St}E. St\o rmer, A survey of noncommutative dynamical entropy, in:
Classification of nuclear C*-algebras. Entropy in operator algebras,
Encyclopaedia Math. Sci., 126, Springer, Berlin, 2002.
\end{thebibliography}
\end{document}